\documentclass[10pt]{article}

\usepackage{amscd,amsmath, amssymb, fancyhdr, mathbbol, dutchcal, url}
\usepackage[matrix,arrow,curve]{xy}
\usepackage[backref=page]{hyperref}

\renewcommand*{\backrefalt}[4]{%
	\ifcase #1 (Not cited.)%
	\or        (Cited on page~#2.)%
	\else      (Cited on pages~#2.)%
	\fi}

\hypersetup{
	colorlinks   = true,
	citecolor    = magenta,
	linkcolor    =blue,
	urlcolor     =magenta	
}

\numberwithin{equation}{section}

\newcommand{\version}{version 1.0,\ \ Feb. 13, 2026}

\def\eqref#1{(\ref{#1})}

\newcommand{\g}{{\mathfrak g}}

\newcommand{\arrow}{{\:\longrightarrow\:}}
\newcommand{\Z}{{\Bbb Z}}
\def\C{{\Bbb C}}
\newcommand{\R}{{\Bbb R}}
\newcommand{\Q}{{\Bbb Q}}

\newcommand{\6}{\partial}

\def\1{\sqrt{-1}\:}
\newcommand{\restrict}[1]{{\left|_{{\phantom{|}\!\!}_{#1}}\right.}}
\newcommand{\cntrct}                
{\hspace{2pt}\raisebox{1pt}{\text{$\lrcorner$}}\hspace{2pt}}

\renewcommand{\tilde}{\widetilde}
\renewcommand{\bar}{\overline}
\renewcommand{\phi}{\varphi}
\renewcommand{\epsilon}{\varepsilon}
\renewcommand{\geq}{\geqslant}

\newcommand{\im}{\operatorname{im}}

\newcommand{\Hom}{\operatorname{Hom}}

\newcommand{\Lie}{\operatorname{Lie}}

\newcounter{Mycounter}[section]
\newcounter{lemma}[section]
\newcounter{claim}[section]
\newcounter{sublemma}[section]
\newcounter{corollary}[section]
\renewcommand{\thecorollary}{{Corollary \thesection.\arabic{corollary}}}
\newcommand{\corollary}{%
    \setcounter{corollary}{\value{Mycounter}}
    \refstepcounter{corollary}
    \stepcounter{Mycounter}
    {\noindent \bf \thecorollary:\ }}

\newcounter{theorem}[section]
\renewcommand{\thetheorem}{{Theorem \thesection.\arabic{theorem}}}
\newcommand{\theorem}{%
    \setcounter{theorem}{\value{Mycounter}}
    \refstepcounter{theorem}
    \stepcounter{Mycounter}
    {\noindent \bf \thetheorem:\ }}

\newcounter{conjecture}[section]
\newcounter{proposition}[section]
\renewcommand{\theproposition}
      {{Proposition \thesection.\arabic{proposition}}}
\newcommand{\proposition}{%
    \setcounter{proposition}{\value{Mycounter}}
    \refstepcounter{proposition}
    \stepcounter{Mycounter}
    {\noindent \bf \theproposition:\ }}

\newcounter{definition}[section]
\renewcommand{\thedefinition}
      {{Definition~\thesection.\arabic{definition}}}
\newcommand{\definition}{%
    \setcounter{definition}{\value{Mycounter}}
    \refstepcounter{definition}
    \stepcounter{Mycounter}
    {\noindent \bf \thedefinition:\ }}

\newcounter{example}[section]
\renewcommand{\theexample}{{Example \thesection.\arabic{example}}}
\newcommand{\example}{%
    \setcounter{example}{\value{Mycounter}}
    \refstepcounter{example}
    \stepcounter{Mycounter}
    {\noindent \bf \theexample:\ }}

\newcounter{remark}[section]
\renewcommand{\theremark}{{Remark \thesection.\arabic{remark}}}
\newcommand{\remark}{%
    \setcounter{remark}{\value{Mycounter}}
    \refstepcounter{remark}
    \stepcounter{Mycounter}
    {\noindent \bf \theremark:\ }}

\newcounter{problem}[section]
\newcounter{question}[section]
\renewcommand{\thequestion}{{Question \thesection.\arabic{question}}}
\newcommand{\question}{%
    \setcounter{question}{\value{Mycounter}}
    \refstepcounter{question}
    \stepcounter{Mycounter}
    {\noindent \bf \thequestion:\ }}

\newcommand{\proof}{\noindent{\bf Proof:\ }}

\newcommand{\pstep}{\noindent{\bf Proof. Step 1:\ }}

\makeatletter

\@addtoreset{equation}{section} \@addtoreset{footnote}{section}

\def\x@arrow{\DOTSB\Relbar}
\def\xlongrightarrowfill@{\arrowfill@\relbar\relbar\longrightarrow}
\newcommand{\xlongrightarrow}[2][]{%
        \ext@arrow 0099\xlongrightarrowfill@{#1}{#2}}
\makeatother

\def\blacksquare{\hbox{\vrule width 5pt height 5pt depth 0pt}}
\def\endproof{\blacksquare}

\begin{document}

\begin{center}
{\LARGE\bf
Dolbeault formality for complex nilmanifolds}\\[4mm]


Tommaso Sferruzza\footnote{The author is holder of a INdAM fellowship and he is partially supported by GNSAGA project ``Progetti di Ricerca 2025 - CUP E53C24001950001".}, Misha Verbitsky,\footnote{Partially supported 
by FAPERJ SEI-260003/000410/2023 and CNPq - Process 310952/2
021-2. 

}

\end{center}

{\small \hspace{0.10\linewidth}
\begin{minipage}[t]{0.85\linewidth}
{\bf Abstract.} 
A quasi-isomorphism of differential graded algebras (DGA)
is a multiplicative map inducing an isomorphism on
cohomology. A DGA is called  formal if it can be connected
by a chain of  quasi-isomorphisms to its cohomology
algebra. We prove that the Dolbeault DGA of a complex nilmanifold is
formal only if it is a torus, and the Dolbeault algebra of
(0,p)-forms is formal if and only if the complex structure
is abelian.
\end{minipage}
}

\tableofcontents


\section{Introduction}

Neisendorfer and Taylor \cite{_NT_} defined a complex
manifold $M$ to be {\bf Dolbeault formal} if the
differential bigraded algebra (shortly, DBA) of bigraded
forms $(\Lambda^{*,*}M,\bar\6)$ is equivalent, as a DBA,
to a DBA with trivial (0,1)-differential.

The equivalence of DBA is understood as a chain of
quasi-isomorphisms (that is, homomorphisms inducing
an isomorphism on the cohomology), preserving the
bigrading.

Nearly fifty years after this definition was given, the
landscape of Dolbeault formal manifolds in the complex
non-K\"ahler setting remains almost entirely unexplored,
safe for few isolated classes of examples.

In this paper, we analyze the existence of Dolbeault
formal manifolds in the class of nilmanifolds endowed with
invariant complex structures.

By a renowned result of Malcev \cite{_Malcev_}, a {\bf
  nilmanifold} is characterized as the compact quotient of
a connected nilpotent Lie group with rational Lie algebra,
by a discrete lattice.

This class of compact manifolds has served as a testing
ground for conjectures and counterexamples in differential
and complex geometry, starting from the work of Thurston
\cite{_Thurston_}.

In the late 80's, Benson and Gordon
\cite{BG_kahler_nilpotent} and, independently, Hasegawa
\cite{_Hasegawa_},  proved that if a nilmanifold admits a
K\"ahler structure than it must be diffeomorphic to a
torus.


The latter argued by proving that the minimal model of a
Lie algebra associated to any nilmanifold is not formal,
unless the Lie algebra is abelian (and, hence, the
nilmanifold is a torus). Indeed, Nomizu's theorem
\cite{_Nomizu_} states that the cohomology of any rational
Lie algebra is isomorphic to the de Rham cohomology of the
associated nilmanifold. On the other hand, compact
K\"ahler manifolds, and more in general, manifolds
satisfying the $dd^c$-lemma, are formal in the sense of
Sullivan \cite{_DGMS_}.

Formal manifolds play a fundamental role in rational homotopy
theory \cite{_Sullivan_,_Felix_Halperin_Thomas_}. Indeed,
for simply connected manifolds (or manifolds with 
nilpotent fundamental group and unipotent action
of the fundamental group on the homotopy groups), the rational
homotopy type of the manifold is a formal consequence of
its cohomology ring.

A strengthening of formality is the notion of {\bf
  geometrically formal metric} \cite{_Kotschick_}, i.e., a
Riemannian metric for which the products of harmonic forms
is still harmonic.  On the other hand, the existence of
non-vanishing Massey products and higher operations
obstructs rational formality.

However, the vanishing of Massey products
does not imply the formality (\cite{_Halperin_Stasheff_}).
This was the motivation for developing the
theory of $A_\infty$-structures, intended
to replace the Massey products in studying 
of rational homotopy (\cite{_Stasheff:A_infty_}). Later,
Babenko and Taimanov suggested a non-commutative
generalizetion of Massey products \cite{_Babenko_Taimanov_}  which contains
a complete information on rational homotopy,
and vanishes if and only if the DGA is formal.
The Massey products for Dolbeault DGA were
computed in \cite{_Cattaneo_Tomassini_},
who gave sufficient conditions for existence
of non-vanishing Massey products.

Several homotopy theories and notions of geometric
formality and Massey products have been introduced
\cite{_NT_,_MS_,_AT_,_TomTor_,_ST3_}.

In 1978, Neisendorfer and Taylor \cite{_NT_}
introduced the Dolbeault homotopy theory,
the homotopy theory for the differential
graded algebras equipped with the Hodge bigrading. 
As in rational homotopy theory, Dolbeault 
formality obstructs the $dd^c$-lemma and, in
particular, the existence of K\"ahler structures on
compact complex manifolds.

To this regard, concerning the Dolbeault formality of
nilmanifolds endowed with invariant complex structure (or,
as we will usually call them, {\bf complex nilmanifolds}),
only partial results have been obtained.

It was announced in 
\cite[Theorem 5]{_Cordero_Fernandez_Gray_Ugarte_2} that nilmanifolds
endowed with a nilpotent complex structures are never
Dolbeault formal, referring to the unpublished preprint
\cite{_Cordero_Fernandez_Gray_Ugarte_3}. In
\cite{_Cordero_Fernandez_Gray_Ugarte_}, the same authors
provide four examples of complex nilmanifolds which are
not Dolbeault formal.

More recently, the first author and Tomassini \cite{_ST3_}
proved the weaker result that complex parallelizable
nilmanifolds are never Dolbeault formal unless they are
complex tori.

A substantial difference with de Rham cohomology of
nilmanifolds arises since there exist no counterpart of
Nomizu's theorem. Indeed, except for certain classes of
invariant complex structures
\cite{_Console_Fino_,_Cordero_Fernandez_Gray_Ugarte_,Rollenske_albanese,_Rollenske:iterated_toric_},
it is not known whether the injection of the Lie algebra
Dolbeault cohomology into the Dolbeault cohomology of the
nilmanifold is an isomorphism.

In this paper, our main theorem finally provides a proof
for any nilmanifold endowed with any invariant complex
structure.

More precisely, let $M=G/\Gamma$ be a complex nilmanifold
and $(\Lambda^{*,*}M,\bar\6)$ its the associated
differential bigraded algebra. Then, the algebras
$(\Lambda^{*}M,\bar\6):=(\text{Tot}^*(\Lambda^{*,*}M),\bar\6)$
and $(\Lambda^{0,*}M,\bar\6)$ are called, respectively,
the {\bf Dolbeault DGA} and {\bf Dolbeault (0,$*$)-DGA} of
$M$. With this in mind, we prove the following 
(\ref{_non-formal_D_Corollary_}).

\hfill 

\theorem
The Dolbeault DGA of a complex nilmanifold $M$ is never
formal, unless $M$ is a torus.

\hfill

The proof relies on the quasi-isomorphism between
the de Rham DGA of a nilmanifold and the Chevalley-Eilenberg DGA
of its Lie algebra, which is due to Nomizu, \cite{_Nomizu_}. 
The Chevalley-Eilenberg DGA, as shown by Hasegawa
(\cite[Theorem 1]{_Hasegawa:non-formal_}), is 
 never formal. Using the 
domination argument due to \cite{_MSZ_},
we are able to use the non-formality
of the Chevalley-Eilenberg DGA to prove that
the Dolbeault DGA of a nilmanifold is non-formal.

Now, let $\g\otimes \C=\g^{1,0}\oplus_\R \g^{0,1}$ be a
complex structure on the rational Lie algebra $\g$
associated to a complex nilmanifold $M=G/\Gamma$. From 
the Serre's duality, triviality of the canonical
bundle and the domination theorem of \cite{_MSZ_}, we obtain the following.

\hfill

\theorem\label{_abelian_*,0_formal_Theorem_}
The Dolbeault (0,$*$)-DGA of a complex nilmanifold is
formal if and only if $\g^{0,1}$ is abelian.

\proof \ref{_formal_then_abelian_Theorem_}. \endproof

\hfill

It would be interesting to study formality of complex
nilmanifolds in the context of pluripotential homotopy
theory \cite{_AT_,_Ste_}. Here, a complex manifold $M$ is
said to be {\bf weakly formal} \cite{_MS_} if the double
complex $(\Lambda^{*,*}M,\6,\bar\6)$ is connected via a
chain of weak equivalences to a bidifferential bigraded
algebra
$(\mathcal{B},\6_{\mathcal{B}},\bar\6_{\mathcal{B}})$ such
that $\6_{\mathcal{B}}\bar\6_{\mathcal{B}}\equiv 0$.

Here, the picture is drastically different from both
rational and Dolbeault homotopy theories. In fact, weak
formality does not obstruct the $dd^c$-lemma nor the
existence of K\"ahler structures
\cite{_ST2_,_Placini_Stelzig_Zoller_}.

Furthermore, on the one hand, the first author and
Tomassini proved that complex nilmanifolds are weakly
formal if and only if they are SKT \cite{_ST_}, in complex
dimension 3. On the other hand, complex parallelizable
solvmanifolds, up to complex dimension 5, are not weakly
formal, \cite{_ST3_} whereas the picture
  is still not clear in higher dimensions. This leaves
open the following question.

\hfill

\question
Which complex nilmanifolds are weakly formal?

\hfill

The paper is organized as follows. In Section 2, we
briefly recall the notions of formality for DG-algebras
and $dd^c$-lemma and the classical proof of formality and
Dolbeault formality of compact K\"ahler manifolds. In
Section 3, recall the definition of nilmanifolds endowed
with complex structures and describe some classical
examples. In Section 4, we recall the definitions of
Chevalley-Eilenberg cohomology and Lie algebra Dolbeault
cohomology and relate them to the de Rham and Dolbeault
cohomologies of nilmanifolds.

In Section 5, we prove the main results of this paper.


\section{Formality for DG-algebras}


\subsection{Quasi-isomorphism of DG-algebras}

\definition
{\bf A graded commutative algebra}, or {\bf
  super-commutative algebra} is a graded algebra 
$A^*= \bigoplus_{i\geq  0} A^i$, with the 
multiplication $A^i\cdot A^j\arrow A^{i+j}$
and graded (super-) commutativity relation
$x\cdot y = (-1)^{\tilde x \tilde y} y\cdot x$,
where $x\in A^{\tilde x}, y \in A^{\tilde y}$.
{\bf A DG-algebra} ({\bf differential
graded algebra}; also abbreviated as DGA) is a graded commutative
algebra equipped with a differential
$d:\; A^i \arrow A^{i+1}$ which satisfies
$d(x\cdot y) = dx \cdot y + (-1)^{\tilde x} x \cdot dy$.
A homomorphism of DG-algebras is called
{\bf a quasi-isomorphism}
if it induces an isomorphism on cohomology.
Two DG-algebras are called {\bf quasi-isomorphic}
if they can be connected by a chain of quasi-isomorphisms.

\hfill

\definition
A DG-algebra is called {\bf  formal}
if it is quasi-isomorphic to its cohomology algebra
(with zero differential). A manifold is called
{\bf formal} if its de Rham algebra is formal.

\subsection{The $dd^c$-lemma}

\definition
Let $M$ be a complex manifold,
and $I:\; TM \arrow TM$ its complex structure
operator. {\bf The twisted differential}
of $M$ is defined as $d^c=IdI^{-1}:\; \Lambda^*(M) \arrow
\Lambda^{*+1}(M)$,
where $I$ acts on 1-forms as an operator dual
to $I:\; TM \arrow TM$, and on the rest of
differential forms multiplicatively.

\hfill

The following theorem, called ``the $dd^c$-lemma'',
was well known. It owes its name and its status
as a foundational result to the paper 
\cite{_DGMS_}, where the formality of K\"ahler
manifolds was first established.

\hfill

\theorem
Let $\eta$ be a form on a compact K\"ahler manifold,
satisfying one of the following conditions.
\begin{enumerate}
\item $\eta$ is an exact $(p,q)$-form. 
\item
$\eta$ is $d$-exact, $d^c$-closed. 
\item $\eta$ is $\bar\6$-exact, and $\6$-closed.
\end{enumerate}
Then $\eta\in \im dd^c$.

\proof \cite{_DGMS_}. \endproof

\hfill

The following corollary, observed in \cite{_DGMS_},
immediately follows from the $dd^c$-lemma.

\hfill

\corollary\label{_q-iso_d^c_Corollary_}
Let $(\Lambda_{d^c}^*(M), d)$ be the algebra of $d^c$-closed forms on 
a compact K\"ahler manifold $M$, considered as a
$DG$-algebra. Then the natural embedding
$\Psi:\; (\Lambda_{d^c}^*(M), d)\arrow (\Lambda^*(M), d)$
is a quasi-isomorphism.

\hfill

\pstep
Any closed $(p,q)$-form is $d^c$-closed.
Since any cohomology class can be represented as a sum
of closed $(p,q)$-forms, it can be represented
by $d^c$-closed form. Therefore $\Psi$ is surjective
on cohomology.

\hfill

{\bf Step 2:}
Suppose that $\Psi$ maps a $d$-closed form
 $\alpha\in\Lambda_{d^c}^*(M)$ to an exact form. 
By $dd^c$-lemma, $\alpha= dd^c\beta$, hence
$\alpha\in d\left(\Lambda_{d^c}^*(M)\right)$.
Therefore, $\Psi$ is injective on cohomology.
\endproof

\hfill

The same way one also proves the following result.

\hfill

\corollary\label{_q-iso_del_Corollary_}
Let $(\Lambda_{\6}^*(M), \bar\6)$ be the algebra of $\6$-closed forms on 
a compact K\"ahler manifold $M$, considered as a
$DG$-algebra. Then the natural embedding
$\Psi:\; (\Lambda_{\6}^*(M), d)\arrow (\Lambda^*(M), \bar\6)$
is a quasi-isomorphism.
\endproof

\subsection{Formality for K\"ahler manifolds}

We are mainly interested in the formality of the Dolbeault
algebra, but the same argument also proves the formality
for all K\"ahler manifolds (a famous result originally
due to \cite{_DGMS_}).

\hfill

\theorem
A compact K\"ahler manifold is formal.

\hfill

\pstep The natural embedding
$(\Lambda_{d^c}^*(M), d)\arrow (\Lambda^*(M), d)$
is a quasi-isomorphism (\ref{_q-iso_d^c_Corollary_}).

\hfill

{\bf  Step 2:} Let 
$\Phi:\; (\Lambda_{d^c}^*(M), d)\arrow (H^*_{d^c}(M),0)$
map a $d^c$-closed form to its cohomology class
in $\frac{\ker d^c}{\im d^c}=: H^*_{d^c}(M)$.
Representing cohomology classes by closed $(p, q)$-forms,
we find that $d\restrict {H^*_{d^c}(M)}=0$. Therefore,
$\Phi$ is a homomorphism of DG-algebras.

It is surjective on cohomology because each class
in $H^*_{d^c}(M)$ can be represented by a $d$, $d^c$-closed
form, and injective because each $d^c$-closed, $d$-exact form
is $dd^c$-exact.
\endproof

\subsection{Dolbeault formality for K\"ahler manifolds}

\definition
Let $(M,I)$ be a compact complex manifold, 
and $\bar\6= \frac{d -\1 d^c}2$ the (0,1)-part of de Rham differential.
The DGA $(\Lambda^*(M), \bar\6)$ is called {\bf Dolbeault algebra,} and
$\bar\6$ {\bf the Dolbeault differential.}
We also consider {\bf  the $(0,*)$-Dolbeault algebra}
$(\Lambda^{0,*}(M), \bar\6):=\left(\bigoplus_{p}\Lambda^{0,p}(M), \bar\6\right)$.

 In \cite{_NT_}, the
  following theorem was proved.

\hfill

\theorem
The Dolbeault algebra 
$(\Lambda^*(M), \bar\6)$ is formal when $(M,I)$ is compact
and K\"ahler.

\hfill

\pstep
Let $(\Lambda^*_\6(M),  \bar\6)$ be the algebra of $\6$-closed forms.
By \ref{_q-iso_del_Corollary_} the natural embedding 
$(\Lambda^*_\6(M),  \bar\6)\to (\Lambda^*(M),  \bar\6)$
is a quasi-isomorphism.

\hfill 

{\bf Step 2:} The map $(\Lambda^*_\6(M),  \bar\6)\to (H^{*}_{\6}(M),0)$
 is a quasi-isomorphism, also by $dd^c$-lemma:
any $\6$-exact, $\bar\6$-closed form is also $\6\bar\6$-exact. 
\endproof


\section{Nilmanifolds}


In this section we define nilmanifolds, list some of their basic
properties, and give a few examples, for the benefit of
the reader. Nothing in this section is original, except
(possibly) the presentation and the choice of examples.

\subsection{Nilmanifolds as homogeneous spaces}

Nilmanifolds are a rich subject deserving much attention,
and complex structures on nilmanifolds provide a rich
source of examples and counterexamples in complex geometry.
In this section we give a brief introduction to nilmanifolds;
for more details and literature see 
\cite{_Vinberg_Gorbatsevich_Shvartsman_,_Corwin_Greenleaf_}.

\hfill

\definition
Let $M$ be a compact manifold equipped
with a transitive action of nilpotent Lie group.
Then $M$ is called {\bf a nilmanifold}.

\hfill

\remark
As shown by A. Maltsev (\cite{_Malcev_}),
all nilmanifolds are obtained as quotient spaces,
$M=G/\Gamma$, where $\Gamma\subset G$ is a cocompact lattice in
a nilpotent Lie group.

\hfill

\theorem\label{_Maltsev_classification_Theorem_}
Let $\g$ be a nilpotent Lie algebra defined over $\Q$,
and $G$ its Lie group.  Then $G$ contains a discrete
subgroup $\Gamma$ such that $G/\Gamma$ is compact, and
$\Gamma= e^{\Gamma_\g}$, where $\Gamma_\g$ is an integer
subalgebra in $\g$. Moreover, $\g \cong \Gamma_\g\otimes_\Q \R$.
Finally, all nilmanifolds are obtained this way.

\proof \cite{_Malcev_}. \endproof

\hfill

\remark Topologically,  all simply connected
nilpotent Lie groups are diffeomorphic to $\R^n$, and all nilmanifolds are 
iterated circle fibrations. This observation (already made
by Maltsev) can be easily deduced from his theorem.

\subsection{Invariant complex structures on nilmanifolds}

\definition
{\bf An integrable complex structure} on a real Lie algebra
$\g$ is a subalgebra $\g^{1,0}\subset \g \otimes_\R \C$
such that $\g^{1,0}\oplus \overline{\g^{1,0}} = \g \otimes_\R \C$

\hfill

\remark  Any such decomposition defines a complex structure $I$
on $\g$ by $I\restrict{\g^{1,0}}=\1$ and $I\restrict{\g^{0,1}}=-\1$.
We extend this complex structure to a left-invariant almost complex
structure on its Lie group $G$. Its integrability is given by 
$[T^{1,0}G, T^{1,0}G]\subset T^{1,0}G$,
which is equivalent to $[\g^{1,0}, \g^{1,0}]\subset \g^{1,0}$.

\hfill

\remark 
In other words, left-invariant complex structures on a
connected real Lie group
 are in bijective correspondence with integrable
complex structures on its Lie algebra.

\hfill

\definition
A {\bf complex nilmanifold} is a nilmanifold
$M=G/\Gamma$ equipped with a complex structure $I$, in such a way
that $G$ has a right-invariant complex structure,
and the projection $G \arrow M$ is holomorphic.
In this case, the complex structure $I$
is called {\bf invariant}. 

\hfill

\remark
Note that not all complex structures on a nilmanifold are invariant;
indeed, a compact torus admits a non-K\"ahler complex
  structure
(\cite{_Blanchard:Recherche_,_Sommese:quaternionic_}), 
which is clearly non-invariant.

\subsection{Parallelizable complex structures}

Consider a complex  nilmanifold $M=G/\Gamma$.
The group $G$ does not need to be a complex Lie group.
When it is complex, and the comple structure on $M$
is induced by the complex structure on this complex
Lie group, we deal with a special class
of complex nilmanifolds called parallelizable.

\hfill

\example An {\bf Iwasawa manifold} is the quotient of
$\begin{pmatrix} 1 & * & * \\
0 & 1 & * \\
0 & 0 & 1
\end{pmatrix}$ 
(group $U_3(\C)$ of upper triangular complex matrices) by
a lattice $\Gamma$. As an example of $\Gamma$ we can take
$U_3(\Z[\1])$.

\hfill

\remark Let $I$ be a bi-invariant complex structure on 
a nilpotent Lie group $G$. Such a complex structure is called
{\bf parallelizable.} For such a complex structure, the
tangent bundle is trivialized by a frame of invariant 
holomorphic vector fields.

\hfill

\remark A complex structure on $\g$ is bi-invariant if and
only if $M:=(G,I)/\Gamma$ is homogeneous under the left
action of $G$. Then the left invariant
vector fields are holomorphic, and this implies that
$[X, Y]=0$ whenever $X\in T^{1,0}M$ and $Y\in T^{0,1}M$.

\hfill

\corollary
Let $I$ be a complex structure on a Lie algebra $\g$.
Then $I$ is bi-invariant if and only if $[\g^{0,1}, \g^{1,0}]=0$.

\hfill

\proof
If this $[\g^{0,1}, \g^{1,0}]=0$, then $\Lie_X Y=0$ is of type (1,0) for any
right-invariant vector fields $X,Y$, with $Y$ of type (1,0).
Then $gYg^{-1}\in \g^{1,0}$ for any $g\in G$, hence
$T^{1,0}G$ is bi-invariant. \endproof

\subsection{Kodaira surface}

\definition
Let $T$, $T'$ be elliptic curves.
{\bf Kodaira surface} $\pi:\; M\arrow T$ is a locally trivial
holomorphic fibration over $T$ with fiber $T'$
and non-trivial Chern class.

\hfill

{\bf A remark on terminology:}
These are {\bf ``primary''}  Kodaira surfaces.
{\bf ``Secondary''} ones are obtained by taking
finite unramified quotients.

\hfill

\remark
The Kodaira surface is diffeomorphic to a quotient
$S^1 \times (G/G_\Z)$ where $G$ is a 3-dimensional 
Heisenberg group, and $G_\Z$ a lattice in $G$.
Therefore, Kodaira surface is a nilmanifold.
This observation is explained in 
\cite{_Hasegawa_}.

\hfill

\remark The Kodaira manifold $M$ is a complex nilmanifold, 
but it is not homogeneous; this would follow e.g. from an
explicit description of its automorphism group given
in \cite{_OV:Principles_}.

\hfill

\remark
Kodaira surface is not K\"ahler. Indeed, 
the cohomology class of $\pi^*(\omega_T)$ vanishes, 
where $\omega_T$ is the K\"ahler form on $T$.
The product of $\omega_T$ and the Kahler form on $M$
(if it exists) is a positive volume form, hence
 cannot be exact.


\section{The Chevalley-Eilenberg complex and its cohomology}


In this section, we define the cohomology of a Lie algebra
(known as the Chevalley-Eilenberg cohomology) and relate
this cohomology algebra with the de Rham and Dobeault
cohomology of a nilmanifold.

\subsection{The cohomology of a Lie algebra}

\proposition\label{_Jacobi_via_CE_diffe_Proposition_}
Let $w\in  \Hom(\Lambda^2 V, V)$.
Consider the map dual to $w$,  $d_w:\; V^* \arrow \Lambda^2 V^*$.
Extend this map to $d_w:\; \Lambda^kV^* \arrow\Lambda^{k+1} V^*$
using the Leibniz rule $d_w(x\wedge y)= d_w(x) \wedge y+
(-1)^{\tilde x} x\wedge d_w y$.  Then $d_w^2=0$
if and only if $w$ defines the Lie algebra structure on
$V$.

\proof \cite{_Chevalley_Eilenberg_}. \endproof

\hfill

\definition
The complex $V^* \stackrel{d_w} \arrow \Lambda^2 V^*  \stackrel{d_w}
\arrow \Lambda^3 V^* \stackrel{d_w}
\arrow ...$ is called {\bf the Chevalley-Eilenberg complex
of the Lie algebra $(V, [\cdot, \cdot]_w)$.}

\hfill

\remark
We identify the Chevalley-Eilenberg complex
of the Lie algebra $\g$ with the complex of left-invariant differential
forms in its Lie group $G$. This defines a natural
homomorphism of DGA $(\Lambda^* \g^*, d_w) \to (\Lambda^*(G/\Gamma), d)$ for
any discrete subgroup $\Gamma\subset G$.

\hfill

\theorem {\bf (Nomizu)}\\
 Let $M=G/\Gamma$ be a nilmanifold.
Then the natural map 
$(\Lambda^* \g^*, d_w) \to (\Lambda^*(G/\Gamma), d)$ 
from its Chevalley-Eilenberg DG-algebra
to its Rham algebra is a quasi-isomorphism. 

{\bf Proof:} \cite{_Nomizu_}. \endproof

\hfill

\corollary
Nilmanifolds, with exception of tori,  are never formal.

\hfill

\proof The Chevalley-Eilenberg DG-algebras 
of finite dimensional nilpotent Lie algebras are never formal,
 \cite[Theorem 1]{_Hasegawa:non-formal_}.
\endproof

\subsection{Dolbeault cohomology of an algebra}

\definition
Let $(\g, I)$ be a complex structure on a nilpotent Lie
algebra, $d_w:\; \Lambda^*(\g^*)\to \Lambda^{*+1}(\g^*)$
its Chevalley-Eilenberg differential, and
$d_w=\6_w+\bar\6_w$ its Hodge decomposition.
The cohomology of the complex $(\Lambda^*_\C(\g^*), \bar\6_w)$
is called {\bf the Dolbeault cohomology of the
Lie algebra}.

\hfill

\remark
In \ref{_Jacobi_via_CE_diffe_Proposition_} we have 
established a correspondence
between the differentials $d_w:\; \Lambda^*(\g^*)\to \Lambda^{*+1}(\g^*)$
and the Lie algebra structures on $\g$. It is natural to ask
which bracket corresponds to $\6_w$ and $\bar\6_w$.

\hfill

\proposition\label{_Dolbeault_bracket_Proposition_}
Let $(\g, I)$ be a complex structure on a nilpotent Lie
algebra, and $\bar\6_w:\; \Lambda^{*,*}(\g^*)\to \Lambda^{*,*+1}(\g^*)$
its Dolbeault differential.  Then $\bar\6_w$ corresponds
to the following bracket $[\cdot, \cdot]_{\bar\6_w}$ on $\g$: given two 
vectors $x, y\in \g^{1,0}$ and two vectors $x',y'\in \g^ {0,1}$,
we have $[x,y]_{\bar\6_w}=0$, $[x,x']_{\bar\6_w}=[x,x']^{1,0}$,
and $[x',y']_{\6_w}=[x',y']$.

\hfill

\proof
The operator $\bar\6_w$ is the (0,1)-part of $d_w$, and $[\cdot, \cdot]_{\bar\6_w}$
is the dual map. In other words, $[\cdot, \cdot]_{\bar\6_w}$
is the Hodge part of the tensor $[\cdot, \cdot]\in \Hom(\Lambda_\C^2 \g, \g\otimes_\R\C)$
which takes a bivector of Hodge type $(1,1)$ to $(1,0)$
and a  bivector of Hodge type $(0,2)$ to $(0,1)$.
\endproof

\hfill



\proposition
Let $(\g,I)$ be a complex structure on a nilpotent Lie algebra, $d_w=\6_w+\bar\6_w$ its Chevalley-Eilenberg differential.  Then $(\g\otimes_\R\C,[\cdot,\cdot]_{\bar\6_w})$ is a complex nilpotent Lie algebra.

\hfill

\proof Since $d_w^2\equiv0$, also $\bar\6_w^2\equiv 0$, which is equivalent to the Jacobi identity for $[\cdot,\cdot]_{\bar\6_w}$. For every $x\in\g\otimes_\R\C$, $\text{ad}(x)$ is nilpotent, so every $\text{ad}_{\bar\6_w}(x)$ is also nilpotent. \endproof

\hfill

\corollary
The DGA $(\Lambda^{*,*}(\g^*), \bar\6_w)$ is formal
if and only $\g$ is abelian. The DGA
$(\Lambda^{0,*}(\g^*), \bar\6_w)$ is formal
if and only $\g^{0,1}$ is abelian.

\proof
Follows immediately from the previous propositions and 
\cite[Theorem 1]{_Hasegawa:non-formal_}.
\endproof

\subsection{Console-Fino theorem}

We start by restating \ref{_Maltsev_classification_Theorem_},
in order to emphasize the rational structure on the
Lie algebra of a nilmanifold.

\hfill

\theorem 
Let $\g$ be a real nilpotent Lie algebra,
and $\Gamma\subset G$ a lattice in its Lie group.
Then $\g$ admits an integer lattice, $\g= \g_\Z \otimes_\Z\R$,
closed under the commutator,
such that $\Gamma$ is commensurable with $e^{\g_\Z}$.
Conversely,  for any such integer subalgebra 
$\g_\Z \subset\g$ such that $\g= \g_\Z \otimes_\Z\R$, the subgroup
$e^{\g_\Z}$ is a lattice in $G$. Finally,  
the rational lattice
$\g_\Q := \g_\Z \otimes_\Z \Q$ is unambiguously determined by
the lattice $\Gamma\subset G$.

\proof \cite{_Malcev_}. \endproof

\hfill

\definition
Let $I$ be an invariant complex structure on a
nilmanifold $G/\Gamma$. We say that $I$ is {\bf  rational}
if it preserves the rational lattice
$\g_\Q \subset \g$ defined by $\Gamma$ as above.

\hfill

\theorem  (Console-Fino)\\
Let $I$ be a rational complex structure on a Lie algebra,
and $(M,I)$ the corresponding complex nilmanifold. Then 
the natural DGA homomorphism $(\Lambda^*(\g), \bar\6_w)
\to (\Lambda^*(M), \bar\6)$ is a quasi-isomorphism.

\proof \cite{_Console_Fino_}. \endproof

\hfill

\remark
In this case, the Dolbeault algebra 
$(\Lambda^*(M), \bar\6)$ is formal if and only if
$(\Lambda^*(\g^*), \bar\6_w)$ is formal.
In the next section, we are going to prove that the Dolbeault algebra
$(\Lambda^*(M), \bar\6)$ is not formal
even when the Console-Fino theorem cannot be applied.


\section{Domination of DG-algebras and formality of Dolbeault DGA}


\subsection{Domination of DG-algebras}
\label{_Domination_Subsection_}

\definition
We say that a DG-algebra $(A^*, d)$ {\bf admits Poincar\'e duality}
if $A^i=0$ for $i <0$ and $i >n$, $H^n(A^*)=\C$,
and the product map $H^i(A^*) \times H^{n-i}(A^*) \to H^n(A^*)=\C$
defines a non-degenerate pairing in cohomology.

\hfill

\definition
Let $\Psi:\; (B^*, d)\to (A^*,d)$ be a morphism of DG-algebras,
such that $(B^*, d)$ admits Poincar\'e duality,
and the map $H^n(B^*) \to H^n(A^*)$ is injective.
Then we say that {\bf  $A^*$ dominates $B^*$} and 
$\Psi$ is called {\bf  a domination}.

\hfill

\example
A proper fibration $f:\; X \to Y$ to an $n$-manifold $Y$ admitting a section
defines a surjection $H_n(X)\to H_n(Y)$, hence
it is injective on cohomology; in this case
$(\Lambda^*(X),d)$ dominates $(\Lambda^*(Y),d)$.

\hfill

\theorem {\bf (Milivojevic, Stelzig, Zoller)}\\
Let $\Psi:\; (B^*, d)\to (A^*,d)$ be a domination of DGA.
Assume that $(B^*, d)$ is not formal.  Then $(A^*, d)$ is not formal.

\proof
\cite[Theorem A]{_MSZ_}. 
\endproof

\subsection{The Dolbeault DGA of a nilmanifold dominates}


Serre's duality immediately implies that 
the Dolbeault algebra $(\Lambda^* M, \bar\6)$
admits a non-degenerate pairing,
taking values in $H^{n,n}_{\bar\6}(M)= \C$.
The same result for the Dolbeault cohomology
of the Lie algebra is also elementary.

\hfill

\proposition\label{_Serre_for_CE_Proposition_}
Let $\g$ be a nilpotent Lie algebra equipped with a 
complex structure. Then $H^{n,n}_{\bar\6}(\Lambda^{*,*}\g^*)=\C$.
Moreover, the multiplication
\[
H^{p,q}_{\bar\6}(\Lambda^{*,*}\g^*)\times H^{n-p,n-q}_{\bar\6}(\Lambda^{*,*}\g^*)
\to H^{n,n}_{\bar\6}(\Lambda^{*,*}\g^*)=\C
\]
defines a non-degenerate pairing.

\proof \cite{_Koszul_}. \endproof

\hfill

Using the terminology introduced in Subsection 
\ref{_Domination_Subsection_},  we say that
these DGAs {\em admit Poincar\'e duality}.

\hfill

\corollary 
Let $(M=G/\Gamma)$ be a complex nilmanifold, and
$\g$ its Lie algebra.  Then the natural
DGA morphism $(\Lambda^*\g^*,\bar\6_w)\to (\Lambda^* M, \bar\6)$
is a domination.

\proof Follows immediately from \ref{_Serre_for_CE_Proposition_}.
\endproof

\hfill

\corollary\label{_non-formal_D_Corollary_}
The Dolbeault DGA of a complex nilmanifold $M$ is never formal,
unless $M$ is a torus.

\proof 
Since $(\Lambda^*\g^*,\bar\6)\to (\Lambda^* M, \bar\6)$
is a domination, and $(\Lambda^*\g^*,\bar\6)$
is not formal, \cite{_MSZ_} imply that $(\Lambda^* M, \bar\6)$
is also not formal.
\endproof

\subsection{Dolbeault formality for $\Lambda^{*,0}(M)$}

It remains to prove \ref{_abelian_*,0_formal_Theorem_}.
A complex structure on a Lie algebra $\g$ is called
{\bf abelian} if the Lie algebra $\g^{1,0}$ is abelian.
There is a great body of research on abelian complex structures
(\cite{_Barberis_Dotti:solvable_,_ABD:classification_,_Cordero_Fernandez_Ugarte:Abelian_}),
which are better understood than the general complex structures
(for instance, the Console-Fino isomorphism holds for
abelian complex structures without extra assumptions,
\cite[Remark 4]{_Console_Fino_}, see 
also \cite{_Cordero_Fernandez_Gray_Ugarte_})

\hfill

\theorem\label{_formal_then_abelian_Theorem_}
Let $(M,I)$ be a nilmanifold with the DGA $(\Lambda^{0,*}(M), \bar\6)$
formal. Then the corresponding complex structure on the Lie algebra is abelian.

\hfill

\proof
As shown in \cite{_Barberis_Dotti_Verbitsky_},
the canonical bundle of a nilmanifold is trivial.
Therefore, Serre's duality implies that 
the DGA $(\Lambda^{0,*}(M), \bar\6)$ 
admits Poincar\'e duality. 
This implies, in particular,
that the natural map 
$(\Lambda^{0,*}(\g), \bar\6)\to (\Lambda^{0,*}(M), \bar\6)$
is a domination. Therefore, formality of
$(\Lambda^{0,*}(M), \bar\6)$ implies formality of
$(\Lambda^{0,*}(\g), \bar\6)$. 
However, $(\Lambda^{0,*}(\g), \bar\6)$ is the
Chevalley-Eilenberg complex of the Lie 
algebra $\g^{0,1}$. By \cite[Theorem 1]{_Hasegawa:non-formal_},
it is formal if and only if it is
abelian. \endproof

\hfill

\hfill

{\bf Acknowledgements:}
We are grateful to Grigory Papayanov for his
invaluable help and consultations.

\hfill

{\small 
    \noindent {\sc Tommaso  Sferruzza\\
    \sc Istituto nazionale di alta matematica "Federico Severi"  (INdAM)\\
    Piazzale Aldo Moro, 5\\
    00185 Roma, Ro - Italy
    \tt \\
    sferruzza@altamatematica.it
}

\hfill

	\noindent \sc Misha Verbitsky\\
		\sc Instituto Nacional de Matem\'atica Pura e
			Aplicada (IMPA) \\ Estrada Dona Castorina, 110\\
			Jardim Bot\^anico, CEP 22460-320\\
			Rio de Janeiro, RJ - Brasil \\
	\tt verbit@impa.br }
\end{document}